\documentclass[12pt]{article}
\usepackage{amssymb}
\usepackage{amsfonts}
\usepackage{amsmath}
\usepackage{amsthm}
\usepackage{eucal}
\usepackage{latexsym}
\usepackage{pstricks,pst-node}
\theoremstyle{plain}

\newtheorem{prop}{Proposition}
\theoremstyle{remark}
\newtheorem{defn}{Definition}
\newtheorem*{remark}{Remark}
\newtheorem{demo}{Example}
\numberwithin{thm}{section}
\numberwithin{prop}{section}
\numberwithin{defn}{section}
\begin{document}
\title{Updown categories}
\author{Michael E. Hoffman\\
\small Dept. of Mathematics, U. S. Naval Academy\\[-0.8ex]
\small Annapolis, MD 21402 USA\\[-0.8ex]
\small and\\[-0.8ex]
\small Max-Planck-Institut f\"ur Mathematik\\[-0.8ex]
\small Vivatsgasse 7, D-53111 Bonn, Germany\\[-0.8ex]
\small \texttt{meh@usna.edu}}
\date{\small February 25, 2004 
\small MR Classifications:  Primary 18B35, 06A07; Secondary 05A17, 05C05}
\maketitle
\def\al{\alpha}
\def\be{\beta}
\def\de{\delta}
\def\ep{\epsilon}
\def\zt{\zeta}
\def\la{\lambda}
\def\si{\sigma}
\def\Ga{\Gamma}
\def\De{\Delta}
\def\Si{\Sigma}
\def\tilde{\widetilde}
\def\<{\langle}
\def\>{\rangle}
\def\up{\uparrow}
\def\dn{\downarrow}
\def\Zm{\mathbf Z/m\mathbf Z}
\def\Zn{\mathbf Z/n\mathbf Z}
\def\Hom{\operatorname{Hom}}
\def\Aut{\operatorname{Aut}}
\def\Ob{\operatorname{Ob}}
\def\id{\operatorname{id}}
\def\src{\operatorname{src}}
\def\trg{\operatorname{trg}}
\def\im{\operatorname{im}}
\def\h{\operatorname{ht}}
\def\A{\mathcal A}
\def\C{\mathcal C}
\def\D{\mathcal D}
\def\K{\mathcal K}
\def\T{\mathcal T}
\def\Y{\mathcal Y}
\def\N{\mathbb N}
\def\F{\mathbb F_q}
\def\k{\Bbbk}
\def\Op{\mathcal O}
\def\Ne{\mathcal N}
\def\H{\mathfrak H}
\def\P{\mathcal P}
\def\U{\mathfrak U}
\def\W{\mathfrak W}
\def\UU{\mathfrak U\mathfrak U}
\def\SU{\mathfrak S\mathfrak U}
\def\<{\langle}
\def\>{\rangle}
\def\sh{
\setlength{\unitlength}{.5 pt}
\begin{picture}(40,20)
\put(10,2){\line(1,0){20}}
\put(10,2){\line(0,1){10}}
\put(20,2){\line(0,1){10}}
\put(30,2){\line(0,1){10}}
\end{picture}}
\begin{abstract}
A poset can be regarded as a category in which there is at most 
one morphism between objects, and such that at most one of 
$\Hom(c,c')$ and $\Hom(c',c)$ is nonempty for $c\ne c'$.  If we keep 
in place the latter axiom but allow for more than one morphism
between objects, we can have a sort of generalized poset in which
there are multiplicities attached to the covering relations, and
possibly nontrivial automorphism groups.  We call such a category
an ``updown category.''  In this paper we give a precise definition
of such categories and develop a theory for them, which incorporates
earlier notions of differential posets and weighted-relation posets.
We also give a detailed account of ten examples, including the
updown categories of integer partitions, integer compositions, planar
rooted trees, and rooted trees.
\end{abstract}
\section{Introduction}
\par
Suppose we have a set $P$ of combinatorial objects, for example
rooted trees, which naturally form a ranked poset (for rooted
trees, the ranking is by the number of non-root vertices).
Each object of $P$ can be constructed in steps from a 
basic object in rank 0, and $p$ covers $q$ in the partial order
if $p$ can be built from $q$ in one step.  (For rooted
trees, the ``basic object'' is the one-vertex tree, and the
building-up process consists of adding a new edge and terminal
vertex to some existing vertex.)  In this situation, there are
naturally two sets of multiplicities on the covering relations
of $P$:  the number of ways to build up $p$ from $q$ is $u(q;p)$,
and the number of ways to tear down $p$ to get $q$ is $d(q;p)$.
(For example, for rooted trees $p,q$ with $p$ covering $q$, 
$u(q;p)$ is the 
number of distinct vertices of $q$ to which a new edge and terminal 
vertex can be added to get $p$, while $d(q;p)$ is the number of 
distinct terminal edges of $p$ that, when removed, leave $q$.)
These multiplicities may be distinct, as in the case of rooted trees
(studied in detail in \cite{H2}), and the difference is related
to the automorphism groups of objects of $P$.
\par
Now a poset can be thought of as a category with at most one
morphism between objects, and at most one of the sets $\Hom(c,c')$
and $\Hom(c',c)$ nonempty when $c\ne c'$.  If we relax the
the first of these conditions, we allow for multiplicities
(if $c\ne c'$) and automorphisms (if $c=c'$.).  In \S2 we
give a precise definition of an updown category, which allows us
to formalize the notions of the previous paragraph.  We also
define a morphism between updown categories, as well as products of
updown categories.  For any updown category $\C$, we define 
``up'' and ``down'' operators $U$ and $D$ on the free
vector space $\k(\Ob\C)$, $\k$ a field of characteristic 0.
\par
In \cite{H1} a theory of universal covers was developed for
weighted-relation posets, i.e., ranked posets in which each
covering relation has a single number $n(x,y)$ assigned to it.
The universal cover of a weighted-relation poset $P$ is the
``unfolding'' of $P$ into a usually much larger weighted-relation
poset $\tilde P$, so that the Hasse diagram of $\tilde P$ is 
a tree and all covering relations of $\tilde P$ have multiplicity 1.  
Although the description of $\tilde P$ had a natural description
in each of the seven examples considered in \cite{H1}, the 
general construction of $\tilde P$ given in \cite[Theorem 3.3]{H1} 
was somewhat unsatisfactory since it involved many arbitrary choices.  
In \S3 we study unilateral updown categories (i.e., updown categories 
with trivial automorphism groups):  these are essentially ``categorified'' 
weighted-relation posets, and the universal-cover construction 
(Proposition 3.3 below) is much more natural in this setting.
\par
In \cite{S1,S2} Stanley developed a theory of differential posets.
Some ideas of this theory were extended to the case of rooted
trees in \cite{H2}.  In \S4 we offer a more general view of
``commutation conditions'' that may be satisfied by the operators
$U$ and $D$ defined in \S2 for any updown category.
\par
The theory developed here is somewhat similar to Fomin's theory
of duality of graded graphs \cite{F1,F2}, but is both
more restrictive and more general:  more restrictive in that the 
functions $u(p;q)$ and $d(p;q)$ must give rise to the same partial 
order, i.e., for any pair $p,q$ we have $u(p;q)=0$ if and only if 
$d(p;q)=0$; and more general in that we consider weaker commutation
conditions than he does.
\par
In \S5 we offer ten examples, which include all those
given in \cite{H1}.  These include the posets of monomials,
necklaces, integer partitions, integer compositions, and 
both planar rooted trees and rooted trees.
\par
The basic idea of this paper was conceived during the academic
year 2002-2003, when the author was partially supported by the
Naval Academy Research Council.  This paper was written
during a stay at the Max-Planck-Institut f\"ur Mathematik during
the following academic year, while the author was on sabbatical 
leave from the Naval Academy.  The author thanks both the Academy 
and the Institut for their support.
\section{Updown categories}
\par
We begin by defining an updown category.
\begin{defn}
An updown category is a small category $\C$ with a rank functor 
$|\cdot|:\C\to\N$ (where $\N$ is the ordered set of natural numbers 
regarded as a category) such that
\begin{itemize}
\item[A1.]
Each level $\C_n=\{p\in\Ob\C:|p|=n\}$ is finite.
\item[A2.]
The zeroth level $\C_0$ consists of a single object $\hat 0$,
and $\Hom(\hat 0,p)$ is nonempty for all objects $p$ of $\C$.
\item[A3.]
For objects $p,p'$ of $\C$, $\Hom(p,p')$ is always finite, and
$\Hom(p,p')=\emptyset$ unless $|p|<|p'|$ or $p=p'$.  In the
latter case, $\Hom(p,p)$ is a group, denoted $\Aut(p)$.
\item[A4.]
Any morphism $p\to p'$, where $|p'|=|p|+k$, factors as a
composition $p=p_0\to p_1\to\dots\to p_k=p'$, where $|p_{i+1}|=
|p_i|+1$;
\item[A5.]
If $|p'|=|p|+1$, the actions of $\Aut(p)$ and $\Aut(p')$ on $\Hom(p,p')$
(by precomposition and postcomposition respectively) are free.
\end{itemize}
\end{defn}
\par
Given an updown category, we can define the multiplicities mentioned
in the introduction as follows.
\begin{defn}
For any two objects $p,p'$ of an updown category $\C$ with $|p'|=|p|+1$, 
define
$$
u(p;p')=\left|\Hom(p,p')/\Aut(p')\right|=
\frac{|\Hom(p,p')|}{|\Aut(p')|}
$$
and
$$
d(p;p')=\left|\Hom(p,p')/\Aut(p)\right|=
\frac{|\Hom(p,p')|}{|\Aut(p)|}.
$$
\end{defn}
It follows immediately from these definitions that
\begin{equation}
u(p;p')|\Aut(p')|=d(p;p')|\Aut(p)| .
\label{udaut}
\end{equation}
\par
We note two extreme cases.  First, suppose $\C_n$ is empty for
all $n>0$.  Then $\C$ is essentially the finite group $\Aut\hat 0$.
Second, suppose that every set $\Hom(p,p')$ has at most one element.
Then $\C$ is a ranked poset with least element $\hat 0$.
\par
Two important special types of updown categories are defined
as follows
\begin{defn} An updown category $\C$ is unilateral if
$\Aut(p)$ is trivial for all $p\in\Ob\C$.  An updown category
$\C$ is simple if $\Hom(c,c')$ has at most one element for all
$c,c'\in\Ob\C$, and the factorization in A4 is unique, i.e.,
for $|c'|>|c|$ any $f\in\Hom(c,c')$ has a unique factorization
into morphisms between adjacent levels.
\end{defn}
Of course simple implies unilateral, but not conversely.
A unilateral updown category is the ``categorification'' of a 
weighted-relation poset in the sense of \cite{H1}; see the next section 
for details.
\par
If $\C$ and $\D$ are updown categories, their product $\C\times\D$ is the 
usual one, i.e. $\Ob(\C\times\D)=\Ob\C\times\Ob\D$ and
$$
\Hom_{\C\times\D}((c,d),(c',d'))=\Hom_{\C}(c,c')\times\Hom_{\D}(d,d') .  
$$
The rank is defined on $\C\times\D$ by $|(c,d)|=|c|+|d|$.  
We have the following result.
\begin{prop}
If $\C$ and $\D$ are updown categories, then so is their product
$\C\times\D$.
\end{prop}
\begin{proof}
Since 
$$
(\C\times\D)_n=\coprod_{i+j=n}\C_i\times\D_j ,
$$
axiom A1 is clear; and evidently $\hat 0=(\hat 0_{\C},\hat 0_{\D})$ 
satisfies A2.  Checking A3 is routine, and
for A4 we can combine factorizations
$$
c=c_0\to c_1\to\dots\to c_k=c'\quad\text{and}\quad
d=d_0\to d_1\to\dots\to d_l=d'
$$
into
$$
(c,d)\to (c_1,d)\to\dots\to (c',d)\to (c',d_1)\to\dots\to (c',d') .
$$
Finally, for A5 note that, e.g.,
$$
\Hom((c,d),(c',d))\cong \Hom(c,c')\times\Aut(d) ,
$$
and the action of $\Aut(c,d)\cong\Aut(c)\times\Aut(d)$ on this set
is free if and only if the action of $\Aut(c)$ on $\Hom(c,c')$ is free.
\end{proof}
\par
We note that the product of two unilateral categories is unilateral,
but the product of simple categories need not be simple:  see Example 2
in \S5 below.
We now define a morphism of updown categories.
\begin{defn}
Let $\C,\D$ be updown categories.  A morphism from $\C$ to $\D$ is
a functor $F:\C\to\D$ with $|F(p)|=|p|$ for all $p\in\Ob\C$, and
such that, 
for any $p,q\in\Ob\C$ with $|q|=|p|+1$, the induced maps
$$
\Aut(p,p)\to \Aut(F(p),F(p)),
$$
$$
\coprod_{\{q': F(q')=F(q)\}} \Hom(p,q')/\Aut(p)\to\Hom(F(p),F(q))/\Aut(F(p)),
$$
and
$$
\coprod_{\{q': F(q')=F(q)\}} \Hom(p,q')/\Aut(q')\to\Hom(F(p),F(q))/\Aut(F(q))
$$
are injective.
\end{defn}
\par
We have the following result.
\begin{prop} Suppose $F:\C\to\D$ is a morphism of updown categories.
If $\D$ is unilateral, then so is $\C$; if $\D$ is simple, then
$\C$ is also simple and $F$ is injective as a function on object
sets.
\end{prop}
\begin{proof}
It follows immediately from Definition 2.4 that
$\C$ must be unilateral when $\D$ is.
Now suppose $\D$ is simple.  Then $\C$ is
unilateral, and it follows from Definition 2.4 that the induced
function
$$
\coprod_{\{q':F(q')=F(q)\}}\Hom(p,q')\to \Hom(F(p),F(q))
$$
is injective when $|q|=|p|+1$:  but $\Hom(F(p),F(q))$ is (at most)
a one-element set, so $F$ must be injective on object sets and $\Hom(p,q)$ 
can have at most one object.  But then unique factorization of morphisms
in $\C$ follows from that in $\D$, so $\C$ is simple.
\end{proof}
\par
One can verify that there is a morphism of updown categories
$\C\to\C\times\D$ given by sending $c\in\Ob\C$ to $(c,\hat 0_{\D})$
whenever $\C$ and $\D$ are updown categories; similarly there is
a morphism $\D\to\C\times\D$.  We denote the $n$-fold cartesian
power of $\C$ by $\C^n$. 
\par
Let $\k$ be a field of characteristic 0, $\k(\Ob\C)$ the free vector space 
on $\Ob\C$ over $\k$.  
We now define ``up'' and ``down'' operators on $\k(\Ob\C)$.
\begin{defn} For an updown category $\C$, let $U,D:\k(\Ob\C)\to\k(\Ob\C)$
be the the linear operators given by
$$
Up=\sum_{|p'|=|p|+1} u(p;p')p'
$$
and
$$
Dp=\begin{cases} \sum_{|p'|=|p|-1} d(p';p) p',& |p|>0,\\
0, & p=\hat 0,\end{cases}
$$
for all $p\in\Ob\C$.
\end{defn}
\begin{prop} The operators $U$ and $D$ are adjoint with respect
to the inner product on $\k(\Ob\C)$ defined by
\begin{equation*}
\<p,p'\>=\begin{cases} |\Aut (p)|,&\text{if $p'=p$},\\
0,&\text{otherwise.}\end{cases}
\end{equation*}
\end{prop}
\begin{proof} Since $\<Up,p'\>=\<p,Dp'\>=0$ unless $|p'|=|p|+1$, it
suffices to consider that case.  Then
$$
\<Up,p'\>=u(p;p')\<p',p'\>=u(p;p')|\Aut(p')|
$$
while 
$$
\<p,Dp'\>=d(p;p')\<p,p\>=d(p;p')|\Aut(p)|,
$$
and the two agree by equation (\ref{udaut}).
\end{proof}
\par
Now we extend the definitions of $u(p;p')$ and $d(p;p')$ to
any pair $p,p'\in\Ob\C$ by setting $u(p;p')=d(p;p')=0$ if 
$\Hom(p,p')=\emptyset$ and
$$
u(p;p')=\frac{\<U^{|p'|-|p|}(p),p'\>}{|\Aut(p')|},
\quad
d(p;p')=\frac{\<U^{|p'|-|p|}(p),p'\>}{|\Aut(p)|}
$$
otherwise.  It is immediate that equation (\ref{udaut}) still holds,
and that
$$
U^k(p)=\sum_{|p'|=|p|+k} u(p;p') p'
$$
and 
$$
D^k(p)=\sum_{|p'|=|p|-k} d(p;p')p'
$$
for any $p\in\Ob\C$.  (However, it is no longer true that
$u(p;q)$ and $d(p;q)$ have any simple relation to $|\Hom(p,q)|$
when $|q|-|p|>1$.)  An important special case of the extended
equation (\ref{udaut}) is
\begin{equation}
\frac{d(\hat 0;p)}{u(\hat 0;p)}=
\frac{|\Aut(p)|}{|\Aut\hat 0|} 
\label{ratio}
\end{equation}
for any object $p$ of $\C$.  If $\Aut\hat 0$ is trivial (as in
all the examples of \S5 below), equation (\ref{ratio}) gives the order
of the automorphism group of $p\in\Ob\C$ as a ratio of multiplicities
(cf. Proposition 2.6 of \cite{H2}).
We also have the following result.
\begin{prop} If $|p|\le k\le |q|$, then
$$
u(p;q)=\sum_{|p'|=k} u(p;p')u(p';q),
$$
and similarly for $u$ replaced by $d$.
\end{prop}
\begin{proof} We have
\begin{align*}
u(p;q) &= \frac{\<U^{|q|-|p|}p,q\>}{|\Aut(q)|}\\
&=\frac1{|\Aut(q)|}\<U^{k-|p|}U^{|q|-k}(p),q\>\\
&=\frac1{|\Aut(q)|}\sum_{|p'|=k}u(p;p')\<U^{k-|p|}p',q\>\\
&=\frac1{|\Aut(q)|}\sum_{|p'|=k}u(p;p')u(p';q)|\Aut(q)|\\
&=\sum_{|p'|=k}u(p;p')u(p';q) ,
\end{align*}
and the proof for $d$ is similar.
\end{proof}
\begin{defn}
For an updown category $\C$, define the induced partial order on $\Ob\C$ by 
setting $p\preceq q$ if and only if $\Hom(p,q)\ne\emptyset$.  
\end{defn}
It follows from Proposition 2.4 that $p\preceq q\iff u(p;q)\ne 0 \iff
d(p;q)\ne 0$.
Henceforth we write $p \lhd q$ if $q$ covers $p$ in the induced partial 
order.
\par
In the unilateral case, equation (\ref{ratio}) is trivial since
$u(p;q)=d(p;q)$ for all $p$ and $q$.  Nevertheless, we have the
following interpretation of the multiplicity in this case.
\begin{prop}
Let $\C$ be a unilateral updown category, $p,q\in\Ob\C$ with
$|q|-|p|=n>0$.  Then $u(p;q)=d(p;q)$ is the number of distinct
chains $(h_1,\dots,h_n)$ so that each $h_i$ is 
a morphism between adjacent levels and $h_nh_{n-1}\cdots h_1$
is a morphism from $p$ to $q$.
\end{prop}
\begin{proof} We use induction on $n$.  The result is immediate
if $n=1$, since in a unilateral updown category
$$
u(p;q)=d(p;q)=|\Hom(p,q)|
$$
when $|q|=|p|+1$.  Now if $N(p,q)$ denotes the number of
chains $(h_1,\dots,h_n)$ as in the statement of the proposition,
it is evident that, for $|q|>|p|+1$, 
$$
N(p,q)=\sum_{r\lhd q} N(p,r)N(r,q) .
$$
But then the inductive step follows from Proposition 2.4.
\end{proof}
\section{Weighted-Relation Posets and Unilateral Updown Categories}
Let $\U$ be the category of updown categories, $\UU$ the full subcategory
of unilateral updown categories.  For a functor $F$  between unilateral 
updown categories $\C$, $\D$, Definition 2.4 reduces to the requirement 
that $F$ preserve rank and that the induced function
\begin{equation}
\coprod_{\{q':F(q')=F(q)\}}\Hom(p,q')\to\Hom(F(p),F(q))
\label{moruu}
\end{equation}
be injective whenever $p,q\in\Ob\C$ with $|q|=|p|+1$.
\par
The notion of a weighted-relation poset was defined in \cite{H1}.
This consists of a ranked poset
$$
P=\bigcup_{n\ge 0} P_n
$$
with a least element $\hat 0\in P_0$, together with nonnegative
integers $n(x,y)$ for each $x,y\in P$ so that $n(x,y)=0$ unless
$x\preceq y$, and
\begin{equation}
n(x,y)=\sum_{|z|=k}n(x,z)n(z,y)
\label{wrp}
\end{equation}
whenever $|x|\le k\le |y|$.
A morphism of weighted-relation posets $P,Q$ is a rank-preserving map 
$f:P\to Q$ such that
\begin{equation}
n(f(t),f(s))\ge \sum_{s'\in f^{-1}(f(s))} n(t,s')
\label{mwrp}
\end{equation}
for any $s,t\in P$ with $|s|=|t|+1$.  
Let $\W$ be the category of weighted-relation posets.
\par
Given an updown category $\C$, it follows from Proposition 2.4
that the weight functions $n(x,y)=u(x;y)$ and $n(x,y)=d(x;y)$
on the poset $\Ob\C$ (with the partial order defined by Definition
2.6) both satisfy equation (\ref{wrp}).  So we have two 
weighted-relation posets based on $\Ob\C$ corresponding
to these two sets of weights.  In fact, we can describe
them functorially.
\par
If $\C$ is an updown category, we can form a unilateral updown 
category $\C^\up$ with $\Ob\C^\up=\Ob\C$, 
and with $\Hom_{\C^\up}(p,p')$ defined as follows.
We declare $\Hom_{\C^\up}(p,p)=\Aut_{\C^\up}(p)$ trivial for all $p$,
and for $|p'|>|p|$ we define $\Hom_{\C^\up}(p,p')$ to be 
the set $\Hom_{\C}(p,p')$ with the equivalence relation generated
by declaring, for any factorization $f=f_nf_{n-1}\cdots f_1$ of
$f\in\Hom_{\C}(p,p')$ into morphisms between adjacent levels,
$f$ equivalent to $\al_nf_n\cdots \al_1f_1$, where $\al_i\in\Aut(\trg f_i)$.
It is routine to check that $\C^\up$ satisfies the axioms of an updown
category, and for $p,p'\in\Ob\C$ with $|p'|=|p|+1$ the multiplicity is
$$
|\Hom_{\C^\up}(p,p')|=\left|\Hom_{\C}(p,p')/\Aut_{\C}(p')\right|=u(p;p') .
$$
Of course $\C^\up$ coincides with $\C$ if $\C$ is unilateral.
\par
Similarly, for any updown category $\C$ there is a unilateral updown
category $\C^\dn$ with $\Ob\C^\dn=\Ob\C$, trivial automorphisms, and
$\Hom_{\C^\dn}(p,p')$ the set $\Hom_{\C}(p,p')$ with the equivalence
relation $f\sim f_n\be_nf_{n-1}\cdots f_1\be_1$
for $f=f_nf_{n-1}\cdots f_1$ a factorization of $f\in\Hom_{\C}(p,p')$ into
morphisms between adjacent levels and $\be_i\in\Aut(\src f_i)$.
Then 
$$
|\Hom_{\C^\dn}(p,p')|=\left|\Hom_{\C}(p,p')/\Aut_{\C}(p)\right|=d(p;p')
$$
for $p,p'\in\Ob\C$ with $|p'|=|p|+1$.
We have the following result.
\begin{prop} There are two functors $\U\to\UU$, taking an updown
category $\C$ to $\C^\up$ and $\C^\dn$ respectively.
\end{prop}
\begin{proof}
We first consider the ``up'' functor.  For a morphism $F:\C\to\D$
of updown categories, we must produce a morphism $F^\up:\C^\up\to\D^\up$
of unilateral updown categories.  But given such a functor $F$, 
Definition 2.4 requires that $F$ preserve rank and that the induced 
function
$$
\coprod_{\{q':F(q')=F(q)\}}\Hom(p,q')/\Aut(p')\to 
\Hom(F(p),F(q))/\Aut(F(q))
$$
be injective for all $p,q\in\Ob\C$ with $|q|=|p|+1$.  This is
exactly the statement that the induced functor $F^\up$ is a
morphism of unilateral updown categories.
The proof for the ``down'' functor is similar.
\end{proof}
\par
Now we pass from unilateral updown categories to weighted-relation
posets.
\begin{prop} There is a functor $Wrp:\UU\to\W$, sending
a unilateral updown category $\C$ to the set $\Ob\C$ with
the partial order of Definition 2.5 and the
weight function $n(x,y)=u(x;y)=d(x;y)$.
\end{prop}
\begin{proof} The only thing to check is the morphisms.  Suppose
$F:\C\to\D$ is a morphism of $\UU$.  Then $F$ defines a function 
on the object sets, and the function (\ref{moruu}) is injective.
Hence
$$
\sum_{\{q':F(q')=F(q)\}}|\Hom(p,q')|\le |\Hom(F(p),F(q))|
$$
and so (since, e.g., $n(p,q')=|\Hom(p,q')|$), inequality
(\ref{mwrp}) holds and $F$ induces a morphism of weighted-relation
posets.
\end{proof}
\par
As defined in \cite{H1},
a morphism $f:P\to Q$ of weighted-relation posets is a covering map 
if $f$ is surjective and the inequality (\ref{mwrp}) is an equality.  
A universal cover $\tilde P$ of $P$ is a cover $\tilde P\to P$ such 
that, if $P'\to P$ is any other cover, then there is a covering 
map $\tilde P\to P'$ so that the composition $\tilde P\to P'\to P$
is the cover $\tilde P\to P$.
In \cite{H1} such a universal cover was constructed for any 
weighted-relation poset $P$.
\par
In fact, the construction of \cite{H1} can be made considerably
simpler and more natural if we work instead with unilateral updown
categories.  We first categorify the definition of covering map.
\begin{defn}
A morphism $\pi:\C'\to\C$ of unilateral updown categories is a 
covering map if $\pi$ is surjective on the object sets and the 
induced function
\begin{equation}
\coprod_{\{q':\pi(q')=\pi(q)\}}\Hom(p,q')\to\Hom(\pi(p),\pi(q))
\label{cover}
\end{equation}
is a bijection for all $p,q\in\Ob\C'$ with $|q|=|p|+1$.
\end{defn}
Then we have the following result.
\begin{prop}
Every unilateral updown category $\C$ has a universal cover $\tilde\C$.
\end{prop}
\begin{proof}
We define $\tilde\C$ to be the category whose level-$n$ objects are
strings $(f_1,f_2,\dots,f_n)$ of morphisms $f_i\in\Hom(c_{i-1},c_i)$,
where $c_i\in\C_i$, and whose morphisms are just inclusions of
strings.  It is straightforward to verify that $\tilde\C$ is a
unilateral updown category (with $\hat 0_{\tilde\C}$ the empty string).  
Define the functor $\pi:\tilde\C\to\C$ by sending the empty
string to $\hat 0\in\Ob\C$, the nonempty string
$(f_1,\dots,f_n)$ of $\tilde\C$ to the target of $f_n$ in $\Ob\C$,
and the inclusion $(f_1,\dots,f_j)\subset (f_1,\dots,f_n)$ to
the morphism $f_nf_{n-1}\cdots f_{j+1}\in\Hom(c_j,c_n)$.
That the induced function (\ref{cover}) is a bijection is a tautology.
\par
Now let $P:\C'\to\C$ be another cover of $\C$:  we must define a
covering map $F:\tilde\C\to\C'$ of unilateral updown categories so 
that $\pi=PF$.  
We proceed by induction on level:  evidently we can get started by
sending the empty string in $\tilde\C_0$ to the element 
$\hat 0$ of $\C'$.  Suppose $F$ is defined through level $n-1$, and
consider a level-$n$ object $(f_1,\dots,f_n)$ of $\tilde\C$.
Let $c_n=\pi(f_1,\dots,f_n)$.  By the induction hypothesis we
have $c_{n-1}'=F(f_1,\dots,f_{n-1})\in\Ob\C'$, and $c_{n-1}=P(c_{n-1}')$
is the target of $f_{n-1}$, hence the source of $f_n$.
Since
$$
P:\coprod_{\{c':p(c')=c_n\}}\Hom(c_{n-1}',c')\to\Hom(c_{n-1},c_n)
$$
is a bijection, there is a unique morphism $g$ of $\C'$ with
$\src(g)=c_{n-1}'$ sent to $f_n:c_{n-1}\to c_n$.  We define $F(f_1,\dots,f_n)$
to be $\trg(g)$, and the image of the inclusion of $(f_1,\dots,f_{n-1})$
in $(f_1,\dots,f_n)$ to be $g$.  This actually defines the functor
$F$ through level $n$, since by the induction hypothesis $F$ assigns
to the inclusion of any proper substring $(f_1,\dots,f_k)$ in
$(f_1,\dots,f_{n-1})$ a morphism $h$ from $F(f_1,\dots,f_k)$ to $c_{n-1}'$
in $\C'$; then $F$ sends the inclusion of $(f_1,\dots,f_k)$ in 
$(f_1,\dots,f_n)$ to $gh$.
\end{proof}
\par
\begin{remark}
If we think of the functor $\pi:\tilde\C\to\C$ as
a function on the object sets, then the number of objects of
$\tilde\C$ that $\pi$ sends to $p\in\Ob\C$ is $u(\hat 0;p)=d(\hat 0;p)$:
this follows from Proposition 2.5.
\end{remark}
\par
The construction of $\C$ in the preceding result is 
functorial:  given a morphism $F:\C\to\D$ of
unilateral updown posets, we have a morphism $\tilde F:\tilde\C\to
\tilde\D$ given by
$$
\tilde F(f_1,f_2,\dots,f_n)=(F(f_1),F(f_2),\dots,F(f_n)).
$$
Also, the updown category $\C$ is evidently simple.
Thus, if $\SU$ is the full subcategory of simple updown
categories in $\U$, then there is a functor $\UU\to\SU$
taking $\C$ to $\tilde\C$.  In fact, we have the following
result.
\begin{prop} The functor $\UU\to\SU$ taking $\C$ to $\tilde\C$
is right adjoint to the inclusion functor $\SU\to\UU$.
\end{prop}
\begin{proof} It suffices to show that
$$
\Hom_{\UU}(\C,\D)\cong \Hom_{\SU}(\C,\tilde\D)
$$
for any simple updown category $\C$ and unilateral updown category
$\D$.  A morphism $F:\C\to\D$ of unilateral
updown categories gives rise to $\tilde F:\tilde\C\to\tilde\D$,
and since $\C$ is simple there is a natural identification 
$\C\cong\tilde\C$, giving us a morphism $\C\to\tilde\D$.  To
go back the other way, just compose with the covering map
$\pi:\tilde\D\to\D$.
\end{proof}
\section{Commutation Conditions}
We shall consider various conditions on the commutator of the
operators $D$ and $U$ introduced in \S2.
In what follows we write $P_i$ for the restriction of the operator
$P$ to level $i$, so $[D,U]_i=D_{i+1}U_i-U_{i-1}D_i$.
\begin{defn} Let $\C$ be an updown category, with operators $D$ and
$U$ as defined above.  We write $I$ for the identity operator on
$\k(\Ob\C)$.
\begin{itemize}
\item[1.]
If $[D,U]=rI$, where $r$ is a scalar, then $\C$ satisfies the absolute
commutation condition (ACC) with constant $r$.
\item[2.]
If $[D,U]_i=(ai+b)I_i$ for constants $a,b$ then $\C$ satisfies the linear 
commutation condition (LCC) with slope $a$.
\item[3.]
If $[D,U]_i=r_iI_i$ for some sequence of scalars $\{r_0,r_1,\dots,\}$,
then $\C$ satisfies the sequential commutation condition (SCC).
\item[4.]
If every element of $\Ob\C$ is an eigenvector for $[D,U]$, then
$\C$ satisfies the weak commutation condition (WCC).
\end{itemize}
\end{defn}
Evidently ACC $\implies$ LCC $\implies$ SCC $\implies$ WCC.  We
can rephrase the preceding definition as follows.  The updown category 
$\C$ satisfies the WCC if there is a function $\ep:\Ob\C\to\k$ such
that $(DU-UD)(c)=\ep(c)c$ for all $c\in\Ob\C$.
Then $\C$ satisfies the ACC if $\ep(c)$ is independent of $c$, the
LCC if $\ep(c)$ is a linear function of $|c|$, and the SCC if $\ep(c)$
is an arbitrary function of $|c|$.
We have the following result about products; cf. Lemma 2.2.3 of
\cite{F2}.
\begin{prop} Let $\C$ and $\D$ be updown categories.
\begin{itemize}
\item[1.] 
If $\C$ satisfies the ACC with constant $r$ and $D$ satisfies the ACC
with constant $s$, then $\C\times\D$ satisfies the ACC with constant
$r+s$.
\item[2.] 
If $\C$ and $\D$ satisfy the LCC with slope $a$, then so does 
$\C\times\D$.
\item[3.] 
If $\C$ and $\D$ satisfy the WCC, then so does $\C\times\D$.
\end{itemize}
\end{prop}
\begin{proof}
Since any element of $\C\times\D$ covering $(c,d)\in\Ob(\C\times\D)$
must have the form $(c',d)$ with $c'$ covering $c$ or $(c,d')$ with
$d'$ covering $d$, we have
$$
U(c,d)=\sum_{|c'|=|c|+1}u(c;c')(c',d)+
\sum_{|d'|=|d|+1}u(d;d')(c,d')=(Uc,d)+(c,Ud),
$$
and similarly for $D$.  If $\C$ and $\D$ satisfy the WCC, we
can calculate that
$$
(DU-UD)(c,d)=((DU-UD)c,d)+(c,(DU-UD)d)=(\ep(c)+\ep(d))(c,d),
$$
from which all three parts follow easily.
\end{proof}
The following result generalizes Proposition 2.4 of \cite{H2}.
\begin{prop} Let $\C$ be an updown category satisfying the WCC, and 
define $\ep:\Ob\C\to\k$ as above.  Then for objects $c_1,c_2$ of $\C$,
$$
\<U(c_1),U(c_2)\>-\<D(c_1),D(c_2)\>=
\begin{cases} 0,&\text{if $c_1\ne c_2$;}\\
\ep(c)|\Aut(c)|,&\text{if $c_1=c_2=c$.}
\end{cases}
$$
\end{prop}
\begin{proof} Calculate using the adjointness of $U$ and $D$.
\end{proof}
\begin{remark} The second alternative of this result can be
written
$$
\sum_{c'\rhd c}u(c;c')^2|\Aut(c')|-\sum_{c''\lhd c}d(c'';c)^2|\Aut(c'')|=
\ep(c)|\Aut(c)|,
$$
or, dividing by $|\Aut(c)|$ and using equation (\ref{udaut}),
\begin{equation}
\sum_{c'\rhd c}u(c;c')d(c;c')-\sum_{c''\lhd c}u(c'';c)d(c'';c)=\ep(c).
\label{esum}
\end{equation}
\end{remark}
In an updown category satisfying the SCC, we can obtain the
kinds of results proved by Stanley for sequentially differential
posets \cite{S2} and by Fomin for $\mathbf r$-graded graphs in \cite{F2}.  
For example, we have the following result by essentially the same proof 
as Theorem 2.3 of \cite{S2} (see also Proposition 2.7 of \cite{H2}).
\begin{prop}
Let $\C$ be an updown category satisfying the SCC, and let $p\in\C_k$.  
Call a word
$w=w_1w_2\cdots w_s$ in $U$ and $D$ a valid $p$-word if the number of
$U$'s minus the number of $D$'s in $w$ is $k$, and, for each
$1\le i\le s$, the number of $D$'s in $w_i\cdots w_s$ does not
exceed the number of $U$'s.  For such a word $w$, let $S=\{i:w_i=D\}$
and
$$
c_i=|\{j:j>i, w_j=U\}|-|\{j:j\ge i, w_j=D\}|,\quad i\in S .
$$
Then for any valid $p$-word $w$,
$$
\<w\hat 0,p\>=d(\hat 0;p)\prod_{i\in S}(r_0+r_1+\dots+r_{c_i}) .
$$
\end{prop}
This result has the following corollary (cf. \cite[Proposition 2.8]{H2} 
and \cite[Theorem 1.5.2]{F2}).
\begin{prop} 
Let $\C$ be an updown category satisfying the SCC, and let
$p\in\C_k$.  Then for nonnegative $a$,
$$
\sum_{|q|=k+a}d(p;q)u(\hat 0;q)=u(\hat 0;p)\prod_{i=0}^{a-1}
(r_0+r_1+\dots+r_{k+i}) .
$$
\end{prop}
\begin{proof} 
Set $w=D^aU^{a+k}$ in the preceding result to get
$$
\<D^aU^{a+k}\hat 0,p\>=d(\hat 0;p)\prod_{i=0}^{a-1}(r_0+r_1+\dots+r_{k+i}) .
$$
Expand out the left-hand side to get
$$
\sum_{|q|=k+a}u(p;q)d(\hat 0;q)=
d(\hat 0;p)\prod_{i=0}^{a-1}(r_0+r_1+\dots+r_{k+i}) .
$$
Now use equation (\ref{udaut}) and divide by $|\Aut p|/|\Aut\hat 0|$ to obtain
the conclusion.
\end{proof}
In the case $p=\hat 0$ the preceding result is
\begin{equation}
\sum_{|q|=a}d(\hat 0;q)u(\hat 0;q)=\prod_{i=0}^{a-1}(r_0+r_1+\dots+r_i) .
\label{wtsum}
\end{equation}
Comparable results in the case where $\C$ merely satisfies the WCC appear
to be much more complicated.
From equation (\ref{esum}) we have
\begin{equation}
\sum_{|q|=1}u(\hat 0;q)d(\hat 0;q)=\ep(\hat 0) ,
\label{ezero}
\end{equation}
generalizing the case $a=1$ of equation (\ref{wtsum}).  
Since $\ep(\hat 0)\hat 0=(DU-UD)\hat 0 =DU\hat 0$, it follows that
$(DU)^n\hat 0=(\ep(\hat 0))^n\hat 0$ for all $n$.  We use this
in proving the following result, which generalizes cases $a=2$ and
$a=3$ of equation (\ref{wtsum}).
\begin{prop} If $\C$ satisfies the WCC, then
\begin{align*}
\sum_{|q|=2}u(\hat 0;q)d(\hat 0;q)&=\ep(\hat 0)^2 +\sum_{|p|=1}u(\hat 0;p)
d(\hat 0;p)\ep(p) \\
\sum_{|t|=3}u(\hat 0;t)d(\hat 0;t)&=\sum_{|p|=1}u(\hat 0;p)d(\hat 0;p)
(\ep(p)+\ep(\hat 0))^2 +\sum_{|q|=2}u(\hat 0;q)d(\hat 0;q)\ep(q) .\\
\end{align*}
\end{prop} 
\begin{proof} For the first part, write $D^2U^2=(DU)^2+D[D,U]U$
and apply it to $\hat 0$:
\begin{align*}
\<D^2U^2\hat 0,\hat 0\> &= \<(DU)^2\hat 0,\hat 0\>+
\sum_{|p|=1}\<u(\hat 0,p)D[D,U]p,\hat 0\> \\
&=(\ep(\hat 0))^2\<\hat 0,\hat 0\> +
\sum_{|p|=1}\<u(\hat 0;p)\ep(p)Dp,\hat 0\> \\
&=(\ep(\hat 0))^2|\Aut\hat 0|+
\sum_{|p|=1}u(\hat 0;p)\ep(p)d(\hat 0;p)|\Aut\hat 0| .
\end{align*}
On the other hand, 
$$
\<D^2U^2\hat 0,\hat 0\>=\<U^2\hat 0,U^2\hat 0\>=
\sum_{|q|=2}\sum_{|p|=2}u(\hat 0;p)u(\hat 0;q)\<p,q\>=
\sum_{|q|=2}u(\hat 0;q)^2|\Aut(q)|
$$
and the first part follows using equation (\ref{udaut}).
\par
To prove the second part, start by applying
$$
D^3U^3=D^2[D,U]U^2+DUD[D,U]U+D[D,U]^2U+D[D,U]UDU+(DU)^3
$$
to $\hat 0$ and proceed similarly, making use of equation (\ref{ezero}).
\end{proof}
\section{Examples}
\par
In this section we present ten examples of updown posets.  Many
of the associated weighted-relation posets appear in the last
section of \cite{H1}.  For the convenience of the reader we
have included a cross-reference to \cite{H1} at the beginning
of each example where it applies.
\par
\begin{demo} 
Let $\C$ be an updown category such that
$\C_1$ consists of a single object $\hat 1$, $\C_n=\emptyset$
for $n\ne 0,1$, and $\Hom(\hat 0,\hat 1)$ has a single element.
The groups $\Aut(\hat 0)$ and $\Aut(\hat 1)$ are trivial since
they act freely on the one-element set $\Hom(\hat 0,\hat 1)$.
Then
$$
(DU-UD)\hat 0= D\hat 1=\hat 0,
$$
and
$$
(DU-UD)\hat 1=-U\hat 0 = -\hat 1,
$$
so $\C$ satisfies the LCC with slope $-2$.  Of course $\C$ is
simple.
\end{demo}
\begin{demo}
(Subsets of a finite set; \cite[Example 1]{H1})
Let $\D=\C^n$, where $\C$ is the updown category of Example 1.
There is an identification of objects of $\D$ with subsets of
$\{1,2,\dots,n\}$:  an $n$-tuple $(c_1,\dots,c_n)$ corresponds to the
set $\{i: c_i=\hat 1\}$.  
The induced partial order is inclusion of sets, and for
$|p'|=|p|+1$ we have 
$$
u(p;p')=d(p;p')=\begin{cases} 1,&\text{if $p$ is a subset of $p'$,}\\
0,&\text{otherwise.}\end{cases}
$$
The category $\D$ is unilateral, but not simple for $n\ge 2$.
In \cite{H1} it is shown that the universal cover $\tilde\D$ is the 
simple updown category whose level-$m$ elements are linearly ordered 
$m$-element subsets of of $\{1,\dots,n\}$, and whose morphisms are
inclusions of initial segments.
From Proposition 4.1,
$\D$ satisfies the LCC with slope $-2$; in fact, it is easy to see
that $(DU-UD)p=(n-2|p|)p$ for any object $p$ of $\D$.  
\end{demo}
\begin{demo}
Let $\C$ be the category with $\C_n=\{[n]\}$, where $[n]=\{1,2,\dots,n\}$ 
(and $[0]=\emptyset$), and let $\Hom([m],[n])$ be the set of injective 
functions from $[m]$ to $[n]$.  Then the axioms are easily seen to hold, 
with $\Aut[n]=\Si_n$ (the symmetric group on $n$ letters).  Since 
$\Hom([n],[n+1])$ has $(n+1)!$ elements,
we have $u([n];[n+1])=1$ and $d([n];[n+1])=n+1$.  More generally, we have
$u([n];[m])=1$ and $d([n];[m])=m!/n!$ for $m\ge n$.  The unilateral
updown poset $\C^\up$ is just the infinite chain $\N$ regarded as a
simple updown category.  On the other hand, a morphism of $\C^\dn$ from
$[n]$ to $[n+1]$ can be thought of as an element of $[n+1]$ (the
element that a representative injective function $[n]\to[n+1]$ misses),
and so a level-$n$ object of $\tilde\C^\dn$ can be identified with a chain
$(i_1,i_2,\dots,i_n)$ of positive integers with $i_j\le j$.
\par
We have $U([n])=[n+1]$ and $D([n])=n[n-1]$, so
$$
(DU-UD)([n])=(n+1)[n]-n[n]=[n],
$$
and thus $\C$ satisfies the ACC with constant 1.
Cf. Example 2.2.1 of \cite{F2}.
\end{demo}
\begin{demo}
(Monomials; \cite[Example 2]{H1})
Let $\D=\C^n$, where $\C$ is the updown category of Example 3.
Objects of $\D$ can be identified with monomials in $n$ commuting
indeterminates $t_1,\dots,t_n$.  The automorphism group of 
$t_1^{i_1}t_2^{i_2}\cdots t_n^{i_n}$ is $\Si_{i_1}\times\Si_{i_2}
\times\dots\times\Si_{i_n}$, and a monomial $u$ precedes a monomial
$v$ in the induced partial order if $u$ is a factor of $v$.
We have
$$
u(1;t_1^{i_1}\cdots t_n^{i_n})=\frac{(i_1+\dots+i_n)!}{i_1!\cdots i_n!}
\quad\text{and}\quad
d(1;t_1^{i_1}\cdots t_n^{i_n})=(i_1+\dots+i_n)! .
$$
The weighted-relation poset $Wrp(\D^\up)$ appears in \cite{H1},
where it is shown that the universal cover $\tilde{\D^\up}$ can be
identified with the simple updown category whose objects are
monomials in $n$ noncommuting indeterminates $T_1,\dots,T_n$, and
whose morphisms are inclusions as left factors.
From Proposition 4.1, $\D$ satisfies the ACC with constant $n$.
Cf. Example 2.2.2 of \cite{F2}. 
\end{demo}
\begin{demo}
(Necklaces; \cite[Example 3]{H1})
For a fixed positive integer $c$, let $\Ne_m$ be the set of 
$m$-bead necklaces with beads of $c$ possible colors.  More
precisely, a level-$m$ object of $\Ne$ is an equivalence
class of functions $f:\Zm\to [c]$, where $f$ is equivalent
to $g$ if there is some $n$ so that $f(a+n)=g(a)$ for all
$a\in\Zm$.  Thus, for $c=2$ the equivalence class 
$$
\{(1,1,2,2),(2,1,1,2),(2,2,1,1),(1,2,2,1)\}\quad
\text{represents the necklace}
\hskip .3in
\psdots[dotstyle=*](.25,.15)(0,.4)
\psdots[dotstyle=o](-.25,.15)(0,-.1)
\psarc(0,.15){.25}{10}{80}
\psarc(0,.15){.25}{100}{170}
\psarc(0,.15){.25}{190}{260}
\psarc(0,.15){.25}{280}{350}
\hskip .2in .
$$
\vskip .1in
\par\noindent
A morphism from the equivalence class of $f$ in $\Ne_m$ to 
the equivalence class of $g$ in $\Ne_n$
is an injective function $h:\Zm\to\Zn$ with $f(a)=gh(a)$
for all $a\in\Zm$, and such that $h$ preserves the cyclic
order, i.e., if we pick representatives of the $h(i)$ in
$\mathbf Z$ with $0\le h(i)\le n-1$, then some cyclic permutation
of $(h(0),h(1),\dots,h(m-1))$ is an increasing sequence.
\par
In \cite{H1} the covering space of the weighted-relation poset 
$Wrp(\Ne^\up)$ is constructed as the set of necklaces with labelled 
beads. It is also shown that for $p\in\Ne_m$,
$$
u(\hat 0,p)=\frac{m!}{|\Aut(p)|} ,
$$
so it follows from equation (\ref{ratio}) that $d(\hat 0,p)=m!$
for all $p\in\Ne_m$.
\par
If $c=1$, then $\Ne_n$ has a single element $p_n$.  Evidently
$U(p_n)=np_{n+1}$ and $D(p_n)=np_{n-1}$, so
$$
(DU-UD)p_n=n(n+1)p_n-n^2p_n=np_n
$$
and $\Ne$ satisfies the LCC with slope 1.
For $c\ge 2$, $\Ne$ does not satisfy the WCC.  For example, when
$c=2$ we have
$$
(DU-UD)(
\hskip .2in
\psdots[dotstyle=*](-.25,.15)(.25,.15)
\psarc(0,.15){.25}{10}{170}
\psarc(0,.15){.25}{190}{350}
\hskip .2in )=
6
\hskip .2in
\psdots[dotstyle=*](-.25,.15)(.25,.15)
\psarc(0,.15){.25}{10}{170}
\psarc(0,.15){.25}{190}{350}
\hskip .2in + 2
\hskip .2in
\psdots[dotstyle=*](-.25,.15)
\psdots[dotstyle=o](.25,.15)
\psarc(0,.15){.25}{10}{170}
\psarc(0,.15){.25}{190}{350}
\hskip .2in .
$$
\end{demo}
\begin{demo} (Integer partitions with unit weights; \cite[Example 5]{H1})
Let $\Y$ be the category with $\Ob\Y$ the set of integer partitions,
i.e., finite sequences $(\la_1,\la_2,\dots,\la_k)$ of positive integers 
with 
$$
\la_1\ge\la_2\ge\dots\ge\la_k .
$$ 
The level of a partition is $|\la|=\la_1+\la_2+\dots+\la_k$; we write
$\ell(\la)$ for the length (number of parts) of $\la$.  The set of morphisms
$\Hom(\la,\mu)$ contains a single element if and only if $\la_i\le\mu_i$
for all $i$.
Then $\Y$ is evidently unilateral but not simple.
The weights $u(\la;\mu)=d(\la;\mu)$ appear
in the ring of symmetric functions:  we have
$$
s_1^ks_{\la}=\sum_{|\mu|=|\la|+k}u(\la;\mu)s_{\mu}
$$
where $s_{\mu}$ is the Schur symmetric function associated with
the partition $\mu$ (for definitions see \cite{Mc}).
In \cite{H1} it is shown that the universal cover $\tilde\Y$ is the 
poset of standard Young tableaux, so $u(\hat 0;\la)=d(\hat 0;\la)$ is
the number of standard Young tableaux of shape $\la$.
\par
That $\Y$ satisfies the ACC with constant 1 is shown in \cite{S1} 
(Corollary 1.4), where $\Y$ is the motivating example of theory of 
differential posets; $\Y$ also appears as Example 1.6.8 of \cite{F2}.
\end{demo}
\begin{demo} Let $\K$ be the category with $\Ob\K$ the set of
integer partitions, and $\Hom(\la,\mu)$ defined as follows.
Let $\la=(\la_1,\dots,\la_n)$ and $\mu=(\mu_1,\dots,\mu_m)$,
always written in decreasing order.
Then a morphism from $\la$ to $\mu$ is an injective
function $f:[n]\to[m]$ such that $\la_i\le \mu_i$ whenever $f(i)=j$.
\par
The partial order induced on $\Ob\K=\Ob\Y$ is the same as that of the
preceding example:  the difference is that we now have nontrivial
automorphism groups and weights on covering relations.
The automorphism group of $\la=(\la_1,\dots,\la_k)$ is the
subgroup of $\Si_k$ consisting of those permutations $\si$
such that $\la_i=\la_j$ whenever $\si(i)=j$.
For a partition $\la$, let $m_k(\la)$ be the number of 
times $k$ occurs in $\la$.  Then for partitions $\la,\mu$ with 
$|\mu|=|\la|+1$, we can describe the weights explicitly as
$$
u(\la;\mu)=\begin{cases} 
1,&\text{if $\mu$ is obtained from $\la$ by
adding a new part of size 1,}\\
m_k(\la),&\text{if $\mu$ is obtained by increasing a 
part of size $k$ in $\la$ to $k+1$,}\\
0,&\text{otherwise.}
\end{cases}
$$
$$
d(\la;\mu)=\begin{cases} 
m_1(\mu),&\text{if $\mu$ is obtained from $\la$ by
adding a new part of size 1,}\\
m_{k+1}(\mu),&\text{if $\mu$ is obtained increasing a
part of size $k$ in $\la$ to size $k+1$,}\\
0,&\text{otherwise.}
\end{cases}
$$
The weights $d(\la;\mu)$ appear implicitly in \cite{Ki} and
explicitly in \cite{K}, where they are referred to as 
``Kingman's branching'':  see especially Figure 4 of \cite{K}.
As noted there, the $d(\la;\mu)$ have an algebraic interpretation
similar to that of the last example:  in the ring of symmetric
functions we have
$$
m_1^km_{\la}=\sum_{|\mu|=|\la|+k}d(\la;\mu)m_{\mu},
$$
where $m_{\la}$ is the monomial symmetric function associated
with $\la$.
\par
The universal cover $\tilde\K^\up$ can be described in terms
of set partitions:  in fact, we can identify elements of 
$\tilde\K_n^\up$ with partitions of $[n]$ so that the covering 
map $\pi:\tilde\K^\up\to\K^\up$ takes a partition $P$ of $[n]$ 
to the integer partition of $n$ given by the block sizes of $P$.
Actually we will identify elements of $\tilde\K_n^\up$ with ordered
partitions $(P_1,\dots,P_k)$ of $[n]$, where 
$$
|P_1|\ge |P_2|\ge\dots\ge|P_1|
$$
and, if $|P_i|=|P_j|$ for $i<j$, $\max P_i<\max P_j$.
We do this by using the construction of Proposition 3.3.
Assign the unique partition of $[1]$ to the morphism from $\hat 0$ to $(1)$,
and suppose inductively that we have assigned an ordered partition
$P=(P_1,\dots,P_k)$ of $[n]$ to the chain $(h_1,\dots,h_n)$ of morphisms 
between adjacent levels of $\K^\up$ from $\hat 0$ to 
$\trg(h_n)=(\la_1,\dots,\la_k)\in\Ob\K_n^\up$ so that $\la_i=|P_i|$.
Let $f\in\Hom_{\K}(\la,\mu)$ be a representative of the equivalence
class $h_{n+1}\in\Hom_{\K^\up}(\la,\mu)$, where $|\mu|=n+1$.
If $\mu$ has length $k+1$, there is a unique element $i\in [k+1]$ not in 
the image of $f$; in this case assign $(P_1,\dots,P_k,\{n+1\})$ to the chain
$(h_1,\dots,h_n,h_{n+1})$.  Otherwise, $\mu$ has length $k$ and there is
a unique $i\in [k]$ such that $\la_i<\mu_{f(i)}$:  in this case,
assign to $(h_1,\dots,h_{n+1})$ the rearrangement of $(P_1',\dots,P_k')$, 
where
$$
P_j'=\begin{cases} P_j\cup\{n+1\},&\text{if $j=i$},\\
P_j,&\text{otherwise,}
\end{cases}
$$
so that $P_i'$ immediately follows $P_m'$, where 
$m=\max\{j<i:|P_j'|\ge |P_i'|\}$.
Evidently the set partition assigned to $(h_1,\dots,h_{n+1})$ projects
to $\mu$ in either case.
\par
Level-$n$ objects of the universal cover $\tilde\K^\dn$ can be described
as sequences $s=(a_1,\dots,a_n)$ such that $m_1(s)\ge m_2(s)\ge\cdots$,
where $m_i(s)$ is the number of occurrences of $i$ in $s$; the covering
map sends $s$ to $(m_1(s),m_2(s),\dots)$.  As in the preceding paragraph,
we can proceed inductively using the construction of Proposition 3.3.
Start by assigning $s=(1)$ to the morphism from $\hat 0$ to $(1)$.
Suppose now we have assigned $s=(a_1,\dots,a_n)$ to a chain of morphisms
$(h_1,\dots,h_n)$ between adjacent levels of $\K^\dn$ from $\hat 0$ to 
$\la=(\la_1,\dots,\la_k)\in\K_n^\dn$ so that $m_i(s)=\la_i$ for 
$1\le i\le k$, and let $h_{n+1}\in\Hom_{\K^\dn}(\la,\mu)$ where $|\mu|=n+1$.  
Now a representative  $f\in\Hom_{\K}(\la,\mu)$ of $h_{n+1}$ must
be ``almost an automorphism'' exchanging parts of equal size with
just one exception:  there is a unique $i\in [\ell(\mu)]$ such
that either $i$ is not in the image of $f$ (in which case $\mu_i=1$),
or else $\la_{f^{-1}(i)}<\mu_i$ (in which case $\mu_i=\la_{f^{-1}(i)}+1$).
Let $S=\{j > i: \la_{f^{-1}(j)}=\mu_i\}$:  note that $S$ is
independent of the choice of $f$.  Now define a permutation $\si$
of $[\ell(\mu)]$ as follows.  If $S=\emptyset$, let $\si$ be the
identity; otherwise, if $S=\{i+1,\dots,l\}$, let $\si(a)=a+1$ for
$i\le a\le l-1$, $\si(l)=i$, and $\si(a)=a$ for $a\notin\{i,\dots,l\}$.
We then assign the sequence $s'=(\si(a_1),\dots,\si(a_k),i)$ to
the chain $(h_1,\dots,h_n,h_{n+1})$.  If $i\notin\im f$, then
$\mu_j=1$ for $j\ge i$ and either $i=\ell(\mu)=k+1$ (if $S$ is empty)
or $l=\ell(\mu)=k+1$ (if it isn't):  either way $\mu$ differs by $\la$
by having 1 inserted in the $i$th position, and $s'$ projects to $\mu$.
If $\mu_i=\la_{f^{-1}(i)}+1$, then we must have $\la_j=\mu_j$ for
$j<i$, and $\mu$ differs from $\la$ in having a part of size $\mu_i-1$
increased by 1.  If $S$ is empty, $\la_{f^{-1}(i)}=\la_i$ and 
$\mu_i=m_i(s')=m_i(s)+1=\la_i+1$.  
Otherwise, $\mu_i=m_i(s')=m_l(s)+1=\la_{f^{-1}(l)}$
and $m_{j+1}(s')=m_j(s)$ for $i\le j\le l-1$.  Either way, $s'$
again projects to $\mu$.
\par
The updown poset $\K$ satisfies the WCC with $\ep(\la)=1+m_1(\la)$.
That is,
\begin{equation}
(DU-UD)(\la)=(1+m_1(\la))\la
\label{eigpart}
\end{equation}
for all partitions $\la$.
To prove this, we introduce the union operation on partitions,
e.g., $(2,1)\cup(3,1,1)=(3,2,1,1,1)$.  If we extend $\cup$ linearly
to $\k(\Ob\K)$, then it is straightforward to show that
\begin{align*}
D(\la\cup\mu) &= D(\la)\cup\mu + \la\cup D(\mu) \\
U(\la\cup\mu) &= U(\la)\cup\mu + \la\cup U(\mu) - \la\cup\mu\cup(1)
\end{align*}
for partitions $\la,\mu$.  An easy calculation then shows
$$
(DU-UD)(\la\cup\mu)=(DU-UD)(\la)\cup\mu+\la\cup(DU-UD)(\mu)-\la\cup\mu,
$$
and since $m_1(\la\cup\mu)=m_1(\la)+m_1(\mu)$, equation (\ref{eigpart})
must hold for $\la\cup\mu$ whenever it holds for $\la$ and $\mu$.
Since (\ref{eigpart}) is easy to show for partitions with one part,
the general result follows by induction on length.
\end{demo}
\begin{demo}
(Integer compositions; \cite[Example 6]{H1})
Let $\C_n$ be the set of integer compositions of $n$, i.e. sequences 
$I=(i_1,\dots,i_p)$ of positive integers with $a_1+\dots+a_m=n$;
as with partitions we write $\ell(I)$ for the length of $I$.
A morphism from $(i_1,\dots,i_p)\in\C_n$ to 
$(j_1,\dots,j_q)\in\C_m$ is an order-preserving injective
function $f:[p]\to[q]$ such that $i_a\le j_{f(a)}$ for all
$a\in [p]$.  Then $\C$ is a unilateral updown category (but
not simple).  The weights $u(I;J)=d(I;J)$ have an algebraic
interpretation similar to that of the preceding two examples,
but here one has to use the ring of quasi-symmetric functions
(for definitions see \cite[Sect. 9.4]{Re}):  if $M_I$ is the monomial
quasi-symmetric function associated with $I$, then
$$
M_1^kM_I=\sum_{|J|=|I|+k} u(I;J)M_J .
$$
\par
The universal cover $\tilde\C$ is constructed in \cite{H1} using 
Cayley permutations as defined in \cite{MF}:  
a Cayley permutation of level $n$ is a length-$n$ sequence 
$s=(a_1,\dots,a_n)$ of positive integers such that any positive 
integer $i<j$ appears in $s$ whenever $j$ does.  
The covering map $\pi:\tilde\C\to\C$ sends a sequence $s$ to the 
composition $(m_1(s),m_2(s),\dots)$.
To relate this to the construction of Proposition 3.3,
we again proceed inductively.  Send the morphism from $\hat 0$
to $(1)$ to the Cayley permutation $(1)$, and suppose 
we have assigned to a chain $(h_1,h_2,\dots,h_n)$ of morphisms
between consecutive levels of $\C$ from $\hat 0$ to $I=(i_1,\dots,i_k)
\in\C_n$ a Cayley permutation $s=(a_1,\dots,a_n)$ that projects to $I$:
note that $\max\{a_1,\dots,a_n\}=k$.
Now let $h_{n+1}\in\Hom(I,J)$ with $J\in\C_{n+1}$.
Then either $\ell(J)=k$ and $h_{n+1}$ is the identity function on $[k]$, 
or $\ell(J)=k+1$.
In the first case, there is exactly one position $q$ where $J$
differs from $I$:  assign to $(h_1,\dots,h_{n+1})$ the Cayley
permutation $s'=(a_1,\dots,a_n,q)$.  Then $m_q(s')=m_q(s)+1=i_q+1$
and $m_i(s')=m_i(s)$ for $i\ne q$, so $s'$ projects to $J$.
In the second case, there is exactly one element $q\in[k+1]$ that 
$h_{n+1}$ misses:  assign $s'=(h_{n+1}(a_1),\dots,h_{n+1}(a_n),q)$ 
to $(h_1,\dots,h_{n+1})$.  Then $\pi(s')=(m_1(s'),m_2(s'),\dots)$ differs
from $I$ only in having an additional 1 inserted in the $q$th place,
and so must be $J$.
\par
The updown category $\C$ satisfies the WCC with 
\begin{equation}
\ep(I)=\ell(I)+2m_1(I)+1 ,
\label{eigcomp}
\end{equation}
where $m_1(I)$ is the number of $1$'s in $I$.  This can be proved 
by induction on length using a method similar to that used for the 
preceding example.  First, it is easy to show that
$$
(DU-UD)(I)=(2m_1(I)+2)(I)
$$
when $\ell(I)=1$.
Now let $I\sqcup J$ be the juxtaposition of compositions $I$ and $J$,
and extend $\sqcup$ linearly to $\k(\Ob\C)$.  Then
\begin{align*}
D(I\sqcup J) &= D(I)\sqcup J+I\sqcup D(J)\\
U(I\sqcup J) &= U(I)\sqcup J+I\sqcup D(U)-I\sqcup (1)\sqcup J
\end{align*}
for compositions $I,J$.  Hence we can calculate that
$$
(DU-UD)(I\sqcup J)=(DU-UD)(I)\sqcup J+I\sqcup (DU-UD)(J)-I\sqcup J .
$$
Since $\ell$ and $m_1$ are additive with respect to the operation
$\sqcup$, it follows that $I\sqcup J$ is an eigenvector of $DU-DU$
with eigenvalue given by equation (\ref{eigcomp}) whenever $I$ and
$J$ are.
\end{demo}
\begin{demo}
(Planar rooted trees; \cite[Example 4]{H1})
Let $\P_n$ consist of functions $f:[2n]\to\{-1,1\}$ so that the 
partial sums $S_i=f(1)+\cdots+f(i)$ have the properties that $S_i\ge 0$
for all $1\le i\le 2n$, and $S_{2n}=0$.  We declare $\Aut(f)$ to
be trivial for all objects $f$ of $\P$, and define a morphism from
$f\in\P_n$ to $g\in\P_{n+1}$ to be an injective, order-preserving function 
$h:[2n]\to[2n+2]$ such that the two values of $[2n+2]$ not in the
image of $h$ are consecutive, and $f(i)=gh(i)$ for $1\le i\le 2n$.
Then $\P$ is a unilateral updown category.  Using the well-known
identification of balanced bracket arrangements with planar
rooted trees, e.g. 
\begin{equation*}
(1,1,-1,1,1,-1,-1,-1)\quad\text{is identified with}\quad
\psline{*-*}(.25,0)(.5,.5)
\psline{*-*}(.5,.5)(.75,0)
\psline{*-*}(.75,0)(.75,-.5)
\hskip .4in ,
\end{equation*}
\vskip .2in
\par\noindent
we can think of $\P$ as the updown category of planar rooted trees;
the level is the count of non-root vertices.  
\par
The weighted-relation poset $Wrp(\P)$ appears as Example 4 in
\cite{H1}, and its universal cover is described as the poset
whose level-$n$ elements are permutations $(a_1,a_2,\dots,a_{2n})$
of the multiset $\{1,1,2,2,\dots,n,n\}$ such that, if $a_i>a_j$
with $i<j$, then there is some $k<j$, $k\ne i$, such that
$a_k=a_i$. (The covering map sends a sequence $s=(a_1,\dots,a_{2n})$
to a sequence of 1's and $-1$'s by sending the first occurrence
of $i$ in $s$ to 1 and the second to $-1$.)  This construction
can be identified with $\tilde\P$ as constructed in Proposition 3.3
in an obvious way.
For example, consider the morphism from $\hat 0=\emptyset$
to $(1,1,-1,1,-1,-1)$ given by the composition $h_3h_2h_1$,
where $h_1=\emptyset$, $h_2=\{(1,1),(2,4)\}$ and 
$h_3=\{(1,1),(2,2),(3,3),(4,6)\}$.  We can code the chain
$(h_1,h_2,h_3)$ by the sequence $(1,2,2,3,3,1)$.
\par
Using the tree language, we can think of $U(t)$ as the sum of all
planar rooted trees obtained by attaching a new edge and
terminal vertex at every possible position of $t$ (a sum with
$2|t|+1$ terms), and $D(t)$ as the sum of all tree obtained
by deleting a terminal edge of $t$.  For example,
$$
U(
\psline{*-*}(.25,0)(.5,.5)
\psline{*-*}(.5,.5)(.75,0)
\hskip .4in ) =
\psline{*-*}(.25,0)(.55,.5)
\psline{*-*}(.55,.5)(.55,0)
\psline{*-*}(.55,.5)(.85,0)
\hskip .4in +
\psline{*-*}(.25,0)(.5,.5)
\psline{*-*}(.5,.5)(.75,0)
\psline{*-*}(.25,0)(.25,-.5)
\hskip .4in +
\psline{*-*}(.25,0)(.55,.5)
\psline{*-*}(.55,.5)(.55,0)
\psline{*-*}(.55,.5)(.85,0)
\hskip .4in +
\psline{*-*}(.25,0)(.5,.5)
\psline{*-*}(.5,.5)(.75,0)
\psline{*-*}(.75,0)(.75,-.5)
\hskip .4in +
\psline{*-*}(.25,0)(.55,.5)
\psline{*-*}(.55,.5)(.55,0)
\psline{*-*}(.55,.5)(.85,0)
$$
\vskip .2in
and
$$
D(
\psline{*-*}(.25,0)(.5,.5)
\psline{*-*}(.5,.5)(.75,0)
\psline{*-*}(.25,0)(.25,-.5)
\hskip .4in ) =
\psline{*-*}(.25,0)(.5,.5)
\psline{*-*}(.5,.5)(.75,0)
\hskip .4in +
\psline{*-*}(.25,0)(.25,.5)
\psline{*-*}(.25,0)(.25,-.5)
\hskip .2in .
$$
\vskip .2in
\par
The updown poset $\P$ satisfies the WCC with 
\begin{equation}
\ep(t)=2|t|+\tau(t)+1 ,
\label{eigprt}
\end{equation}
where $\tau(t)$ is the number of terminal vertices of $t$.  This
can be proved by a method similar to that of the preceding two
examples, but here we need two operations:  a binary operation
$\vee$ and a unary operation $B_+$.  The binary operation 
$\vee:\P_n\times\P_m\to\P_{n+m}$
can be described as juxtaposition of balanced bracket arrangements,
or equivalently as joining two planar rooted trees at the root:
$$
\psline{*-*}(.25,0)(.5,.5)
\psline{*-*}(.5,.5)(.75,0)
\psline{*-*}(.75,0)(.75,-.5)
\hskip .5in \vee
\hskip .2in
\psline{*-*}(0,0)(0,.5)
\hskip .3in =
\hskip .3in
\psline{*-*}(0,0)(.3,.5)
\psline{*-*}(.3,.5)(.3,0)
\psline{*-*}(.3,0)(.3,-.5)
\psline{*-*}(.3,.5)(.6,0)
\hskip .3in .
$$
\vskip .2in
\par\noindent
The unary operation $B_+:\P_n\to\P_{n+1}$ encloses a balanced
bracket operation in an outer pair of delimiters, or equivalently
adds a new root vertex at the top of a planar rooted tree:
$$
B_+(
\psline{*-*}(.25,0)(.5,.5)
\psline{*-*}(.5,.5)(.75,0)
\hskip .4in ) =
\psline{*-*}(.25,-.5)(.5,0)
\psline{*-*}(.5,0)(.75,-.5)
\psline{*-*}(.5,0)(.5,.5)
\hskip .3in .
$$
\vskip .2in
Now it is straightforward to show that
\begin{align*}
D(t_1\vee t_2) &= D(t_1)\vee t_2 + t_1\vee D(t_2)\\
U(t_1\vee t_2) &= U(t_1)\vee t_2 + t_1\vee U(t_2) - t_1\vee
\psline{*-*}(.25,.35)(.25,-.15)
\hskip .2in \vee t_2
\end{align*}
for any two planar rooted trees $t_1,t_2$, and that
\begin{align*}
D(B_+(t)) &= B_+(D(t))\\
U(B_+(t)) &= B_+(U(t)) + 
\psline{*-*}(.25,.35)(.25,-.15)
\hskip .15in \vee\: U(t) +
U(t)\,\vee
\psline{*-*}(.25,.35)(.25,-.15)
\end{align*}
for any planar rooted tree $t$. Using the first pair of
these equations, we can calculate that
$$
(DU-UD)(t_1\vee t_2)=(DU-UD)(t_1)\vee t_2 + t_1\vee (DU-UD)(t_2)
-t_1\vee t_2
$$
for any planar rooted trees $t_1,t_2$:  since $\tau(t_1\vee t_2)=
\tau(t_1)+\tau(t_2)$, it follows that $t_1\vee t_2$
is an eigenvector of $DU-UD$ satisfying equation (\ref{eigprt}) whenever
$t_1$ and $t_2$ are.  Similarly, the second pair of equations gives
$$
(DU-UD)B_+(t)=B_+((DU-UD)(t))+2B_+(t)
$$
for any planar rooted tree $t$.  Since $\tau(B_+(t))=\tau(t)$, it follows
that $B_+(t)$ is an eigenvector of $DU-UD$ satisfying (\ref{eigprt})
when $t$ is.  Now any planar rooted tree $t$ with $|t|>0$ can be
written as either $t_1\vee t_2$ or $B_+(t_1)$, so we can prove the 
result by induction on $|t|$.
\end{demo}
\begin{demo}
(Rooted trees; \cite[Example 7]{H1})
Let $\T_n$ consist of partially ordered sets $P$ such that
(1) $P$ has $n+1$ elements; (2) $P$ has a greatest element; and
(3) for any $v\in P$, the set of elements of $P$ exceeding $v$
forms a chain.
The Hasse diagram of such a poset $P$ is a tree with the greatest 
element (the root vertex) at the top.
A morphism of $\T$ from $P\in\T_m$ to $Q\in\T_n$ is an injective
order-preserving function $f:P\to Q$ that sends the root of $P$ to the
root of $Q$, and which preserves covering relations (i.e., if
$v\lhd w$ in the partial order on $P$, then $f(v)\lhd f(w)$ in
the partial order on $Q$).  Then $\T$ is an updown category.
\par
The updown category $\T$ was studied extensively in \cite{H2},
though without using the categorical language.  To see that the
construction of the preceding paragraph gives the same multiplicities
as in \cite{H2}, consider a morphism from $P\in\T_n$ to $Q\in\T_{n+1}$.  
Any such morphism misses only some terminal vertex $v\in Q$, so we
can think of it as identifying $P$ with $Q-\{v\}$.  Elements of
$$
\Hom(P,Q)/\Aut(Q)
$$
amount to different choices for the parent of $v$ in $Q$, i.e.,
different choices for terminal vertices of $P$ to which a new
edge and vertex can be attached to form $Q$:  this is $n(P;Q)$
as defined in \cite{H2}.  On the other hand, elements of 
$$
\Hom(P,Q)/\Aut(P)
$$
amount to different choices of $v$, and thus to different
choices for an edge of $Q$ that when cut leaves $P$:  this is
$m(P;Q)$ as defined in \cite{H2}.  
\par
The operators $U$ and $D$ on $\k(\Ob\T)$ appear in \S2 of \cite{H2} as
$\mathfrak N$ and $\mathfrak P$ respectively.  As is proved there 
(Proposition 2.2), $\T$ satisfies the LCC with $\ep(t)=|t|+1$
(Note that the grading in \cite{H2} differs by 1 from the one
used here.)
\par
In \cite{H1} the weighted-relation poset $Wrp(\T^\up)$ is discussed,
and it is shown that the universal cover $\tilde\T^\up$ can be described 
as permutations of $[n]$.  
Finding a simple description of the objects of $\tilde\T^\dn$
appears to be a harder problem.
\end{demo}


\begin{thebibliography}{99}
\bibitem{F1}
S. V. Fomin, The generalized Robinson-Schensted-Knuth correspondence
(Russian), \emph{Zap. Nauchn. Sem. Leningrad. Otel Mat. Inst. Steklov.
(LOMI)} {\bf 155} (1986), Differentsialnaya Geometriya, Gruppy Li i
Mekh. VIII, 156-175; translation in \emph{J. Soviet Math.} {\bf 41}
(1988), 979-991.
\bibitem{F2} 
S. Fomin, Duality of graded graphs, \emph{J. Algebraic Combin.} {\bf 3}
(1994), 357-404.
\bibitem{H1}
M. E. Hoffman, An analogue of covering space theory for ranked posets,
\emph{Electron. J. Combin.} {\bf 8} (2001), res. art. 32.
\bibitem{H2}
M. E. Hoffman, Combinatorics of rooted trees and Hopf algebras,
\emph{Trans. Amer. Math. Soc.} {\bf 355} (2003), 3795-3811.
\bibitem{K}
S. Kerov, The boundary of Young lattice and random Young tableaux,
\emph{Formal Power Series and Algebraic Combinatorics (New Brunswick,
NJ 1994)}, DIMACS Series in Discrete Mathematics and Theoretical 
Computer Science 24, American Mathematical Society, Providence, RI,
1996, pp. 133-158.
\bibitem{Ki}
J. F. C. Kingman, Random partitions in population genetics,
\emph{Proc. Roy. Soc. London Ser. A} {\bf 361} (1978), 1-20.
\bibitem{Mc}
I. G. Macdonald, \emph{Symmetric Functions and Hall Polynomials},
2nd ed., Oxford University Press, New York, 1995.
\bibitem{MF}
M. Mor and A. S. Fraenkel, Cayley permutations, \emph{Disc. Math.}
{\bf 48} (1984), 101-112.
\bibitem{Re}
C. Reutenauer, \emph{Free Lie Algebras}, Oxford University Press,
New York, 1993.
\bibitem{S1}
R. P. Stanley, Differential posets, \emph{J. Amer. Math. Soc.} {\bf 1}
(1988), 919-961.
\bibitem{S2}
R. P. Stanley, Variations on differential posets, \emph{Invariant Theory
and Tableaux (Minneapolis, MN 1988)}, IMA Volumes in Mathematics and its
Applications 19, Springer-Verlag, New York, 1990, pp. 145-165.
\end{thebibliography}
\end{document}